\documentclass[11pt]{article}
\usepackage{amssymb}
\usepackage{amsthm}
\usepackage{amsmath}
\usepackage{graphicx}

\def\fpd#1#2{{\displaystyle\frac{\partial #1}{\partial #2}}}

\def\vf#1{{\displaystyle\frac{\partial }{\partial #1}}}

\def\onehalf{{\textstyle\frac12}}

\newcommand{\R}{\mathbb{R}}

\theoremstyle{plain}
\newtheorem*{thm}{Theorem}

\author{M.\ Crampin and  T.\ Mestdag\\
Department of Mathematics, Ghent University\\
Krijgslaan 281, B--9000 Gent, Belgium\\[0.5cm] {\em This paper is dedicated to Professor Lajos Tam\'{a}ssy}\\{\em in celebration of his 90th birthday.}}

\title{A class of Finsler surfaces whose geodesics are circles}

\begin{document}
\date{}
\maketitle

\begin{abstract}
We determine all Finsler metrics of Randers type for which the Riemannian part is a scalar multiple of the Euclidean metric, on an open subset of the Euclidean plane, whose geodesics are circles. We show that the Riemannian part must be of constant Gaussian curvature, and that for every such Riemannian metric there is a class of Randers metrics satisfying the condition, determined up to the addition of a total derivative, depending on a single parameter. As one of several applications we exhibit a Finsler metric whose geodesics are the oriented horocycles in the Poincar\'{e} disk, in each of the two possible orientations.
\end{abstract}

\subsubsection*{MSC}
53C60

\subsubsection*{Keywords}
Finsler metric, geodesic, circle, horocycle.

\section{Introduction}
When, in two dimensions, does a Finsler metric have geodesics which are all circles? Here `circle' is used in the mildly generalized sense in which straight lines are considered as circles of infinite radius. Apart from this qualification the word has its usual meaning, so there is an underlying assumption of an Euclidean structure:\  and in fact we will work in Euclidean coordinates in the Euclidean plane.

There are several reasons, in addition to pure curiosity, why one might be interested in this question. It may be considered as a generalized version, or in Tabachnikov's words \cite{Tab} as a `magnetic analog', of Hilbert's fourth problem. Or one might be interested in an inverse problem, that is, in the question of whether some particular path space consisting of circles is Finsler metrizable. Our motivation originated in a problem of the latter type, the relevant path space being the set of horocycles of a fixed orientation in the Poincar\'{e} disk. (It is well known that any two-dimensional path space is locally Finsler metrizable. We therefore emphasise that the question here is a global one:\ whether there is a Finsler function --- positive and strictly convex --- defined over the whole of the interior of the unit disk, whose geodesics are horocycles.)

We do not attempt to answer the question in all generality, but deal with the special class of Finsler metrics consisting of Randers metrics for which the Riemannian part is conformal to the Euclidean metric. The obvious two-dimensional Riemannian spaces whose geodesics are circles belong to this class. So do the metrics whose geodesics are circles of fixed radius discussed by Shen in \cite{Shen}, which we have elsewhere \cite{mult} discussed at some length under the name of Shen's circles, and which form the standard case of a magnetic flow as discussed in \cite{Tab}. Now the transformations of the plane which map circles to circles are the fractional-linear or M\"{o}bius transformations; these are also conformal for the Euclidean metric, and so this class of Finsler functions is mapped to itself by the M\"{o}bius transformations. This makes the class a natural one to investigate. We solve the problem completely for this class of metrics. Our results may be summarized as follows.

\begin{thm} Let $M$ be an open subset of the Euclidean plane, with Euclidean coordinates $(x,y)$, and let $(u,v)$ be the corresponding fibre coordinates on the slit tangent bundle $T^\circ\! M$. In order that a Finsler metric of Randers type whose Riemannian part is conformal to the Euclidean metric should have only circles as geodesics the Riemannian part must be of constant Gaussian curvature. Up to a M\"obius transformation, and up to multiplication by a positive constant and the addition of an arbitrary total derivative, such a Finsler metric belongs to one of following three one-parameter families of Randers metrics for which the Gaussian curvature of the Riemannian part is respectively $0$, $+1$ and $-1$:
\begin{align*}
F_\tau(x,y,u,v)&=\sqrt{u^2+v^2}+\tau(yu-xv)
\\
F_\tau(x,y,u,v)&= \frac{\sqrt{u^2+v^2}+ \tau(yu-xv)}{2(1+(x^2+y^2))} 
\\
F_\tau(x,y,u,v)&= \frac{\sqrt{u^2+v^2}+ \tau(yu-xv)}{2(1-(x^2+y^2))}. 
\end{align*}
\end{thm}

We discover, among the third of these classes of Finsler metrics, one for the horocycles in the Poincar\'{e} disk in each of the two possible orientations. This discovery resolves a minor controversy. It is known unequivocally (see \cite{Recht}) that neither of the horocycle flows on a compact surface of constant negative curvature (obtained as the quotient space of the Poincar\'{e} disk by the action of a suitable discrete group of isometries) is geodesible:\ that is to say, no Riemannian metric for such a surface can be found whose geodesic flow has the horocycles (of either orientation) as its geodesic paths. Two questions arise:\ does this result apply to the horocycle flows for the whole disk, and does it rule out Finsler metrizability? A paper by Sullivan \cite{Sull} seems to suggest that it does apply to the whole disk; moreover, the arguments adduced there apparently rule out Finsler metrizability too. Our result shows that in fact each horocycle path space on the disk is Finsler metrizable.

Also among our Finsler metrics are ones for curves with constant geodesic curvature with respect to the Poincar\'{e} metric, which again are circles:\ these are the base integral curves of what Arnol'd \cite{Arnold}  calls `cyclical flows'. It was known that cyclical flows are  isomorphic to geodesic flows, but so far as we know no explicit metric has been published until now.

We should make it clear that we do not claim that our examples exhaust the possibilities for Finsler metrics whose geodesics are circles. Tabachnikov, in \cite{Tab}, describes the whole class of metrics whose geodesics are circles all of the same radius. \'{A}lvarez Paiva and Berck, in \cite{AP}, find all reversible Finsler metrics on the 2-sphere whose geodesics are circles. In each case the methods used are different from ours. Our solutions do of course overlap with those of these other authors, in some instances of the first and second cases listed in the Theorem above, respectively.

In Section 2 we obtain the Randers spaces whose geodesics are circles, thus proving the Theorem. In Section 3 we discuss horocycles. In Section 4 we compute the geodesic curvature of the geodesics of the Randers metrics of interest. In Section 5 we consider cyclical flows.

We have specified the Finsler metrics under consideration as those Randers metrics for which the Riemannian part is conformal to the Euclidean metric. Now in two dimensions every Riemannian metric is locally conformally flat, so it may be argued that that part of the specification is superfluous:\ that we are dealing with Randers metrics in two dimensions for which the geodesics are circles when expressed in terms of coordinates for which the Riemannian part is a scalar multiple of the Euclidean metric. It seems to us however to be more natural geometrically to state the problem in the original form, that is, to start with the Euclidean plane with Euclidean coordinates as  described above, and insist that the Finsler functions in question are of Randers type with Riemannian part which is a scalar multiple of the Euclidean metric in those coordinates.

\section{Proof of the Theorem}
Let $M$ be an open subset of the Euclidean plane, with Euclidean coordinates $(x,y)$, and let $(u,v)$ be the corresponding fibre coordinates on $T^\circ\! M$. The vertical vector fields
\[
u\vf{u}+v\vf{v}\quad\mbox{and}\quad -v\vf{u}+u\vf{v}
\]
form a basis for vertical vector fields. So any spray may be written
\[
\Gamma=u\vf{x}+v\vf{y}+\lambda\left(-v\vf{u}+u\vf{v}\right) \mod(\Delta)
\]
where $\lambda$ is a function on $T^\circ\! M$, positively homogeneous of degree 1 in $(u,v)$. Here
\[
\Delta=u\vf{u}+v\vf{v},
\]
the Liouville field. We wish to determine conditions under which the base integral curves of $\Gamma$ are circles, considered as oriented paths. That is to say, we may count as equivalent sprays which differ by a multiple of $\Delta$, or in other words work projectively. In the first instance, therefore, we may simply ignore the $\Delta$ component. The base integral curves of $\Gamma$ are then the solutions of
\[
\ddot{x}=-\lambda(x,y,\dot{x},\dot{y})\dot{y},\quad \ddot{y}=\lambda(x,y,\dot{x},\dot{y})\dot{x}.
\]
Along a solution curve $\dot{x}^2+\dot{y}^2$ is constant, so we may assume parametrization with Euclidean arc-length. But then $\lambda$ is just the curvature of the solution curve:\ to be more precise, if $s\mapsto (x(s),y(s))$ is a solution curve, $\lambda(x(s),y(s),\dot{x}(s),\dot{y}(s))$ is its curvature at the point of parameter value $s$, the unit tangent to the curve at that point being $(\dot{x}(s),\dot{y}(s))$ and the unit normals $\pm(-\dot{y}(s),\dot{x}(s))$. Now a circle is just a curve of constant curvature:\ so the solution curve is a circle if and only if $\lambda$ is constant along it. Thus the necessary and sufficient condition for all the base integral curves of $\Gamma$ to be circles is that $\Gamma(\lambda)=0$, or that $\lambda$ satisfies
\[
u\fpd{\lambda}{x}+v\fpd{\lambda}{y}+\lambda\left(-v\fpd{\lambda}{u}+u\fpd{\lambda}{v}\right)=0.
\]

We wish next to use this criterion to determine when the geodesics of a Finsler function are circles. It is perhaps too much to expect to be able to answer this question in all generality, but as we have argued in the Introduction there is an obvious class of Finsler functions which are of interest, appropriate to the problem, and easy to analyse:\ it consists of Finsler functions obtained from a conformally-Euclidean Riemannian metric by a Randers change:\ that is, Finsler functions $F$ of the form
\[
F(x,y,u,v)=e^{\phi(x,y)}\sqrt{u^2+v^2}+a(x,y)u+b(x,y)v,
\]
where we have written the conformal factor, which must be everywhere positive, as $e^{\phi(x,y)}$.

Since $F$ is positively homogeneous of degree 1 it has a single Euler-Lagrange equation, which is
\[
\frac{(v\dot{u}-u\dot{v})}{(u^2+v^2)^{3/2}}+
\left(\fpd{\phi}{y}u-\fpd{\phi}{x}v\right)\frac{1}{\sqrt{u^2+v^2}}
+e^{-\phi}\left(\fpd{a}{y}-\fpd{b}{x}\right)=0.
\]
This equation is satisfied by $\dot{u}=-\lambda v$, $\dot{v}=\lambda u$ where
\[
\lambda=e^{-\phi}\left(\fpd{a}{y}-\fpd{b}{x}\right)\sqrt{u^2+v^2}+\fpd{\phi}{y}u-\fpd{\phi}{x}v.
\]
Let us denote the coefficient of $\sqrt{u^2+v^2}$ by $\mu$. Then the condition on $\lambda$ for the geodesics to be circles becomes
\begin{eqnarray*}
\lefteqn{\left(u\left(\fpd{\mu}{x}-\mu\fpd{\phi}{x}\right)
+v\left(\fpd{\mu}{y}-\mu\fpd{\phi}{y}\right)\right)\sqrt{u^2+v^2}}\\
&&\mbox{}
+\left(\frac{\partial^2\phi}{\partial x\partial y}-\fpd{\phi}{x}\fpd{\phi}{y}\right)(u^2-v^2)\\
&&\mbox{}
-\left(\frac{\partial^2\phi}{\partial x^2}-\frac{\partial^2\phi}{\partial y^2}-\left(\fpd{\phi}{x}\right)^2
+\left(\fpd{\phi}{y}\right)^2\right)uv=0.
\end{eqnarray*}
Thus from the quadratic terms $\phi$ must satisfy the two equations
\[
\frac{\partial^2\phi}{\partial x\partial y}=\fpd{\phi}{x}\fpd{\phi}{y}\quad\mbox{and}\quad
\frac{\partial^2\phi}{\partial x^2}-\frac{\partial^2\phi}{\partial y^2}=\left(\fpd{\phi}{x}\right)^2
-\left(\fpd{\phi}{y}\right)^2,
\]
while from the term in $\sqrt{u^2+v^2}$ it follows that $\mu e^{-\phi}$ is constant.

The equations that $\phi$ must satisfy turn out to be simplest when expressed in terms of $e^{-\phi}=f$:\ they become
\[
\frac{\partial^2f}{\partial x\partial y}=0\quad
\mbox{and}\quad
\frac{\partial^2f}{\partial x^2}=\frac{\partial^2f}{\partial y^2}.
\]
It follows that all third order partial derivatives of $f$ vanish, so $f$ is at most quadratic in $x$ and $y$, with highest-order term a constant multiple of $x^2+y^2$. Then the Riemannian part is of constant Gaussian curvature; in fact if $f(x,y)=l(x^2+y^2)+mx+ny+p$ the  curvature  is $4lp-(m^2+n^2)$.

The condition on $a$ and $b$ is that
\[
f^2\left(\fpd{a}{y}-\fpd{b}{x}\right)=k,
\]
a constant. Note that if $(a,b)$ and $(a',b')$ both satisfy this condition (with the same $k$) then
\[
\vf{y}(a'-a)=\vf{x}(b'-b),
\]
so that
\[
a'=a+\fpd{\psi}{x},\quad b'=b+\fpd{\psi}{y}
\]
for some function $\psi$ on $M$:\ but $F$ is only determined up to a total derivative anyway. So any one solution $(a,b)$ is enough. Indeed, if $(a,b)$ is a solution with $k=1$ then $(\tau a,\tau b)$ is a solution with $k=\tau$. In fact with $f(x,y)=l(x^2+y^2)+mx+ny+p$, provided that $m$, $n$ and $p$ are not all zero,
\[
a(x,y)=(\alpha x+\beta y+\gamma)/f(x,y),\quad
b(x,y)=(-\beta x+\alpha y+\delta)/f(x,y)
\]
will do for any constants $\alpha$, $\beta$, $\gamma$ and $\delta$ such that
\[
 m\alpha+n\beta-2l\gamma=0,\quad -n\alpha+m\beta+2l\delta=0,\quad
2p\beta-n\gamma+m\delta=1.
\]

By using M\"{o}bius transformations and the fact that $F$ and $cF$ have the same geodesics for nonzero constant $c$ we can obtain canonical forms, up to the addition of total derivatives, for the Finsler metrics corresponding to Gaussian curvature 0, $+1$, $-1$, as follows.  Note first that since $F$ and $cF$ must both be positive we must take $c$ to be positive; so we may assume without loss of generality that the leading term in $f$ is $(x^2+y^2)$,  $-(x^2+y^2)$ or zero (i.e.\ that $f$ is linear). In either of the first two cases we may eliminate the linear terms by a change of origin, and by scaling reduce the quadratic expression to one of $1+(x^2+y^2)$, $1-(x^2+y^2)$, $(x^2+y^2)-1$ and $(x^2+y^2)$ (the conformal factor cannot be everywhere negative). The latter two cases are equivalent to those with  $f=1-(x^2+y^2)$ and $f=1$ respectively, by inversion in the unit circle. Any linear $f$ may be reduced to $f=y$ (by a rotation followed by a scale followed by a translation), and any constant $f$ to 1. Moreover, it is well-known that there exists a M\"{o}bius transformation in the hyperbolic plane which transforms the half-plane model (with $f=y$)  into the Poincar\'e disk model (with $f=1-(x^2+y^2)$). So we may assume that $f=1$, $f=1+(x^2+y^2)$, or $f=1-(x^2+y^2)$. The first has Gaussian curvature $0$ of course. Finally we replace the latter two by $f=2(1+(x^2+y^2))$ and $f=2(1-(x^2+y^2))$, so as to obtain Gaussian curvature $+1$ and $-1$ respectively. It is then easily checked that $\beta=1$, $\alpha=\gamma=\delta=0$ satisfies the conditions for $a$ and $b$ above, with $k=1$, in each case.  We obtain the following classes of canonical forms
\begin{align*}
F_\tau(x,y,u,v)&=\sqrt{u^2+v^2}+\tau(yu-xv)\qquad\mbox{(curvature $0$)}\\
F_\tau(x,y,u,v)&= \frac{\sqrt{u^2+v^2}+ \tau(yu-xv)}{2(1+(x^2+y^2))}\qquad\mbox{(curvature $+1$)} \\
F_\tau(x,y,u,v)&= \frac{\sqrt{u^2+v^2}+ \tau(yu-xv)}{2(1-(x^2+y^2))}\qquad\mbox{(curvature $-1$)}
\end{align*}
as stated in the Theorem. Each class depends on a single parameter $\tau$.

We should also determine maximum domains on which $F_\tau$ is positive and strongly convex:\ actually positive is enough, since according to \cite{Bao} (Section 11.1), for a Randers metric positivity implies strong convexity. Notice that
\[
|yu-xv|\leq\sqrt{x^2+y^2}\sqrt{u^2+v^2}
\]
by the Cauchy-Schwarz inequality. It follows that the first two Finsler functions above are positive on the domain $x^2+y^2 < \tau^{-2}$. The geodesics of the first class are Shen's circles and we refer to \cite{mult} and \cite{Shen} for further discussions on this metric. As these references reveal, the situation is complicated by the fact that we may change, and possibly enlarge, domains by adding on total derivatives. For the second case, the Finsler function with $\tau=0$ corresponds (up to a factor) to the metric of the sphere after stereographic projection. If $\tau\neq 0$, i.e.\ when the Finsler function is non-reversible, the metrics cannot be extended to the whole sphere. This is indeed to be expected in the light of the following result, which can be derived from \cite{asym}:\ any Finsler metric (not necessarily of Randers type) defined on the whole sphere, all of whose geodesics are circles, is the sum of a reversible metric and a total derivative. We shall discuss the third class of metrics in more detail in the sections that follow.


\section{Horocycles}
As we pointed out earlier, one interesting example consists of the horocycles.

In the Poincar\'{e} disk model of the hyperbolic plane, the horocycles are represented by circles tangent to the boundary circle. So let us take the closed unit disk $D$ centred at the origin in $\R^2$, with interior $D_0$ and boundary the unit circle $C$, and consider all circles whose centres lie in $D_0$ and which touch $C$. Actually, the set of such circles does not define a path space:\ one has first to choose an orientation (otherwise there are two circles through any given point in $D_0$ with given tangent vector). We choose the anticlockwise orientation.

Consider a circle through a point $(x_0,y_0)\in D_0$, with (non-zero) tangent vector $(u_0,v_0)$ there. Its centre $(\xi,\eta)$ lies on the normal, so $x_0-\xi=\rho v_0$,  $y_0-\eta=-\rho u_0$,
with $\rho>0$ to get the orientation right. For this circle to touch $C$ its radius $r$ must satisfy $\sqrt{\xi^2+\eta^2}+r=1$. The equation of the horocycle through $(x_0,y_0)$ with tangent $(u_0,v_0)$ is thus $(x-\xi)^2+(y-\eta)^2=r^2$
where $\xi=x_0-\rho v_0$, $\eta=y_0+\rho u_0$, and
\[
\rho\sqrt{u_0^2+v_0^2}=r=1-\sqrt{(x_0-\rho v_0)^2+(y_0+\rho u_0)^2},
\]
from which
\[
\rho= \frac{1-(x_0^2+y_0^2)}{2\left(\sqrt{u_0^2+v_0^2}+y_0u_0-x_0v_0\right)}.
\]
Since $(x_0,y_0)\in D_0$ the argument of the previous section involving the Cauchy-Schwarz inequality implies that $\rho$ is well-defined and positive. 

We shall show that
\[
F(x,y,u,v)= \frac{\sqrt{u^2+v^2}+yu-xv}{2(1-(x^2+y^2))}
\]
is a Finsler function on $T^\circ\! D_0$ whose geodesics are the horocycles.  Notice that $4F(x_0,y_0,u_0,v_0)=\rho^{-1}$; thus $F$ is positive on $T^\circ\! D_0$, so it is a Finsler function there. It belongs to the class of Finsler functions discussed in the previous section whose geodesics are circles, and whose Riemannian part has Gaussian curvature $-1$; we shall show that its geodesics are the anticlockwise-oriented horocycles.

In this case, somewhat remarkably, $\lambda=4F$. The geodesics are thus the solutions of
\[
\ddot{x}=-4F(x,y,\dot{x},\dot{y})\dot{y},\quad
\ddot{y}=4F(x,y,\dot{x},\dot{y})\dot{x},
\]
with $\dot{F}=0$. Thus $\dot{x}-u_0=-(y-y_0)/\rho$, $\dot{y}-v_0=(x-x_0)/\rho$,
with $\rho=(4 F(x_0,y_0,u_0,v_0))^{-1}$; or in other words
\[
\rho\dot{x}=-y+(y_0+\rho u_0)=-(y-\eta),\quad \rho\dot{y}=x-(x_0-\rho v_0)=x-\xi.
\]
Thus
\[
(x-\xi)^2+(y-\eta)^2=\rho^2(\dot{x}^2+\dot{y}^2)=\rho^2,
\]
since without loss of generality we may take $\dot{x}^2+\dot{y}^2=1$. This base integral path is a circle with centre $(\xi,\eta)$ and radius $\rho$. Finally, we show that $(\xi,\eta)\in D_0$, and $\rho+\sqrt{\xi^2+\eta^2}=1$; in fact the first follows from the second. It follows easily from the definitions of $\rho$ (with $u_0^2+v_0^2=1$) and $\xi$ and $\eta$ that $\xi^2+\eta^2=(1-\rho)^2$. But  by the Cauchy-Schwarz inequality $1-\rho>0$, and therefore $\rho+\sqrt{\xi^2+\eta^2}=1$ as required.

In principle our calculations determine only the projective class of $\Gamma$:\ but $\Gamma(F)=0$ as we showed earlier, so $\Gamma$ is in fact the geodesic spray of $F$.

The geodesics of the Finsler function
\[
\bar{F}(x,y,u,v)= \frac{\sqrt{u^2+v^2}-yu+xv}{2(1-(x^2+y^2))}
\]
are the horocycles with clockwise orientation.

The Finsler metric $F$ for the horocycles has constant flag curvature; but somewhat unexpectedly its flag curvature is 1 rather than $-1$. Since $F$ is a Randers metric we can use the Zermelo navigation data to carry out the calculation. The method is described  by Bao and Robles in \cite{BaoRob}; see also \cite{Robles}. We follow their notation. The Zermelo navigation data for $F$ consist of a Riemannian metric $h$ and a vector field $W$ on $D_0$. Using the formulas given in Section 3.1.2 of \cite{BaoRob} we find that
\[
h = (1-(x^2+y^2))^{-1}\left((1-y^2)dx^2+2xy\, dxdy+(1-x^2)dy^2\right),
\]
while
\[
W=y\vf{x}-x\vf{y}.
\]
Now $h$ is the metric of a hemisphere, obtained by parallel projection of
the upper unit hemisphere in $\R^3$, given by $z=\sqrt{1-(x^2+y^2)}$, along the $z$-axis onto $D_0$. Moreover, $W$ is the generator of rotations about the origin in $\R^2$, so it is an infinitesimal isometry of $h$. It then follows from Theorem~10 of \cite{BaoRob} that $F$ has constant flag curvature equal to the Gaussian curvature of $h$, which is 1.

\section{Geodesic curvature}

Consider again the Finsler function
\[
F(x,y,u,v)=e^{\phi(x,y)}\sqrt{u^2+v^2}+a(x,y)u+b(x,y)v.
\]
One of the conditions we derived earlier for the geodesics of $F$ to be circles is that
\[
e^{-2\phi}\left(\fpd{a}{y}-\fpd{b}{x}\right)
\]
is constant. Let us denote the expression displayed above by $\kappa$. We shall show that for any $F$ of this form and for any geodesic path of $F$, $\kappa$ is its geodesic curvature with respect to the Riemannian part of $F$.

For any $F$ of this form let us denote its Riemannian part $e^{\phi}\sqrt{u^2+v^2}$ by $G$. Any spray $\Gamma$ of $F$ is given by
\[
\Gamma=u\vf{x}+v\vf{y}+\lambda\left(-v\vf{u}+u\vf{v}\right)+\nu\Delta
\]
for some $\nu$, where
\[
\lambda=e^{\phi}\kappa\sqrt{u^2+v^2}+\fpd{\phi}{y}u-\fpd{\phi}{x}v
=\kappa G+\lambda_0
\]
say. Furthermore, let us choose $\nu$ such that $\Gamma(G)=0$, so that the base integral curves of $\Gamma$ are parametrized proportionally to the Riemannian arc-length. This requires just that
\[
\nu=-\left(u\fpd{\phi}{x}+v\fpd{\phi}{y}\right),
\]
from which it is clear that, with this choice of parametrization,
\[
\Gamma=\Gamma_0+\kappa G\left(-v\vf{u}+u\vf{v}\right)
\]
where $\Gamma_0$ is the geodesic spray of the Riemannian part of $F$. Let us write $\nabla/ds$ for the operation of covariant differentiation along a curve parametrized by $s$ with respect to the Levi-Civita connection of the Riemannian metric. Then any base integral curve of the spray $\Gamma$ (with the chosen value of $\nu$), parametrized with Riemannian arc-length so that $G=1$ along it, satisfies
\[
\frac{\nabla^2x}{ds^2}=-\kappa\frac{\nabla y}{ds},\quad
\frac{\nabla^2y}{ds^2}=\kappa\frac{\nabla x}{ds}.
\]
It is clear that $(-\nabla y/ds,\nabla x/ds)$ is a unit normal field along the curve. Thus $\kappa$ is the geodesic curvature of the base integral curve with respect to the Riemannian part of $F$.

Since $\kappa$ is a function on the base $M$ the geodesics of the Finsler functions $F$ under consideration have the interesting property that their geodesic curvatures depend on position only: at each point $(x,y)\in M$ all geodesics of $F$ through $(x,y)$ have the same geodesic curvature $\kappa(x,y)$ there. The flow on $T^\circ\! M$ of the corresponding spray $\Gamma$ (normalized so that $\Gamma(G)=0$) is accordingly called by Arnol'd, in \cite{Arnold}, an isotropic flow. Remark that when $\kappa$ is constant the base integral curves of $\Gamma$ are `Riemannian' circles (in the sense of e.g.\ \cite{Nomizu}), for the Riemannian metric given by $G$.

Consider next the 1-form $\theta=-(adx+bdy)$. We have
\[
d\theta=\left(\fpd{a}{y}-\fpd{b}{x}\right)dx\wedge dy=e^{2\phi}\kappa dx\wedge dy=\kappa \omega,
\]
where $\omega$ is the area 2-form of the Riemannian part of $F$. Let $c$ be a piecewise-smooth closed curve in $M$ which bounds a region $R$:\ then
\[
\oint F(\dot{c})=\oint G(\dot{c})+\int_R\kappa\omega.
\]
Suppose that $\kappa$ is constant, as it is for those Finsler functions $F$ whose geodesics are circles; suppose that it is nonzero; then by scaling $a$ and $b$ we may take $\kappa=1$. Then  the area of $R$ is just $\int_R\kappa\omega=\oint_c \theta$. This observation is relevant to the so-called isoperimetric problem for the Riemannian manifold with metric $G$, that is, the problem of minimising the length of the perimeter, $\oint G(\dot{c})$, of the region while fixing its area. By the method of the Lagrange multiplier a curve $c$ is a solution curve for this isoperimetric problem if it satisfies the Euler-Lagrange equations for the Lagrangian
\[
F_{-\tau}(x,y,u,v)=e^{\phi(x,y)}\sqrt{u^2+v^2}-\tau(a(x,y)u+b(x,y)v)
\]
($\tau$ being the multiplier). That is, solutions of the isoperimetric problem (if any exist) are closed curves which are geodesics of $F_{-\tau}$; the value of $\tau$ is determined by the requirement that the area enclosed by $c$ (as measured with $\omega$) takes the specified value. In particular, we obtain the known results that the solutions to the isoperimetric problem for arbitrary $G$ have constant geodesic curvature, and that the solutions for surfaces of constant curvature are circles. In fact we can generalize the problem slightly, by using $F$ to measure the perimeter rather than $G$, since (for constant $\kappa$) the difference between the two perimeters $\oint G(\dot{c})$ and $\oint F(\dot{c})$ is just a constant multiple of the area enclosed, which is assumed fixed. This discussion complements a solution by Bryant to a problem posted recently by \'{A}lvarez Paiva in an online discussion forum \cite{online}.

\section{Cyclical flows}
We consider now
\[
F_\tau(x,y,u,v)= \frac{\sqrt{u^2+v^2}+ \tau(yu-xv)}{2(1-(x^2+y^2))}
\]
where $\tau$ is a constant. The quantity $f^2\mu$ for this function takes the (constant) value $\tau$:\ its geodesics are therefore circles. If we set $\tau=0$ we get the geodesics of the Poincar\'e disk, if we set $\tau=+1$ or $\tau=-1$ we get horocycles; the configurations for other values of $\tau$ interpolate between, and extend beyond, these. Now in order that $F_\tau>0$ we must have $x^2+y^2<\min(1,\tau^{-2})$. Thus for $-1\leq\tau\leq 1$, $F_\tau$ is a well-defined Finsler metric (positive and strictly convex) over the whole of the open unit disk $D_0$. But for $|\tau|>1$ $F_\tau$ is a Finsler metric only on a smaller domain. It is however what we have elsewhere \cite{mult} called a pseudo-Finsler function; in particular, its geodesics are still well-defined, and are circles.

The possibilities are illustrated in the following figures. Each of the plots show geodesics through the same initial point and with the same initial direction, but for different choices for the parameter $\tau$. The black circles represent either horocycles ($\tau = \pm 1$) or Riemannian geodesics ($\tau = 0$). The grey circles outside the horocycles correspond to Finsler geodesics with parameter values $|\tau|<1$, while those inside the horocycles correspond to $|\tau|>1$.

\begin{figure}[h]
\hspace*{0.5cm}\begin{tabular}{cc}
\begin{minipage}{6cm}
\includegraphics[scale=0.5]{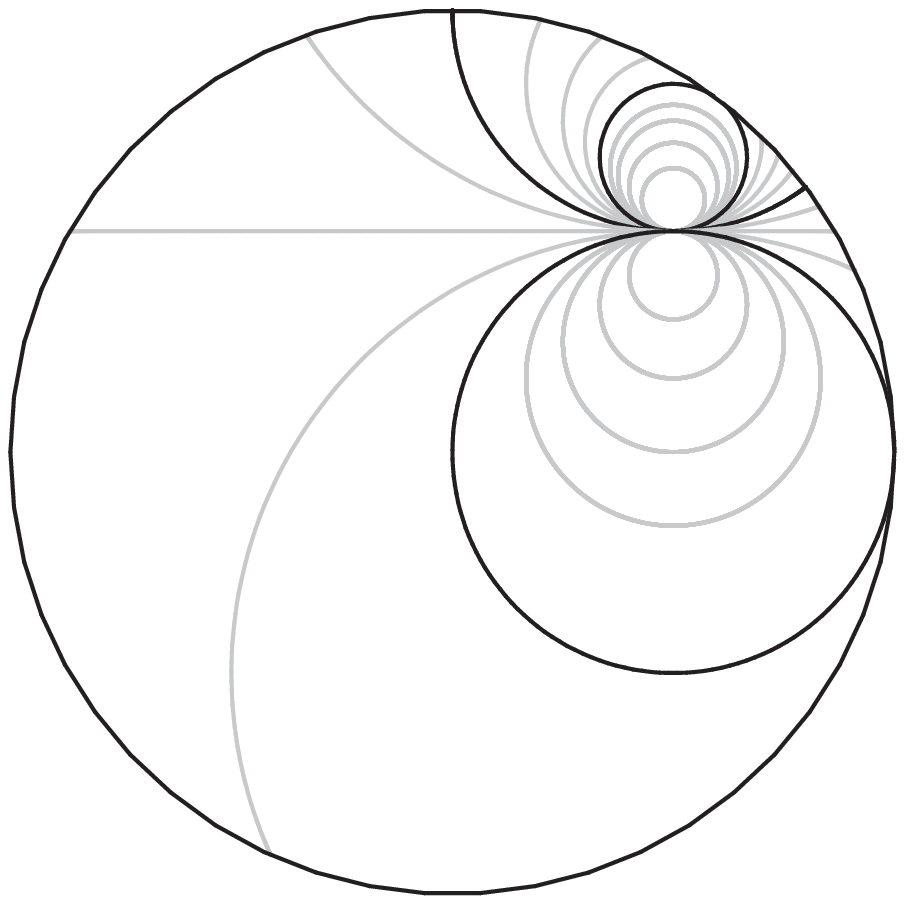}
\end{minipage}
&
\begin{minipage}{7cm}
\includegraphics[scale=0.5]{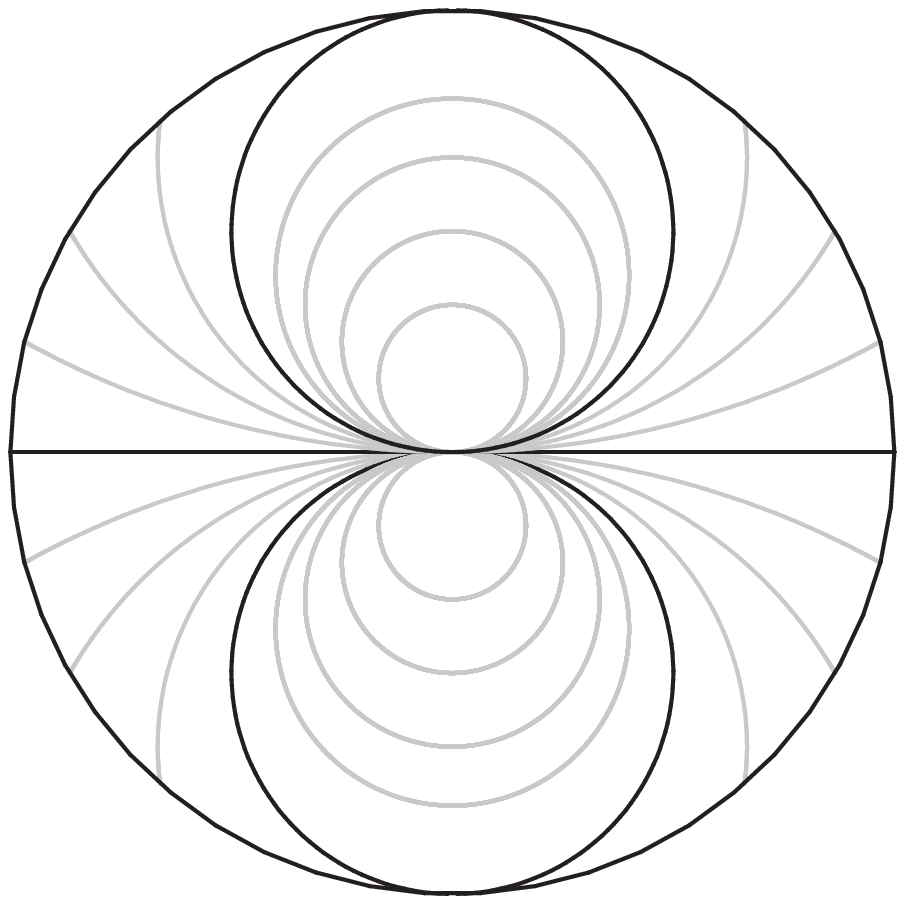}
\end{minipage}
\end{tabular}
\end{figure}

By the calculation in the previous section the geodesic curvature of the integral curves of the geodesic spray of the Finsler function $F_\tau$, viewed as curves on the Poincar\'e disk, is just $\tau$. We are therefore in the case of what Arnol'd calls cyclical flows (\cite{Arnold}; see also \cite{Pat}). Theorem 2 of Arnol'd's paper states that if $\tau^2<1$ the corresponding flow is `isomorphic to a geodesic flow':\  we have exhibited a Finsler function of which it is the geodesic flow. Let us be a little more precise. 
The Finsler function $F_\tau$ determines a {\em Riemannian\/} metric on $T^\circ\! M$, its so-called Sasaki lift. The flow of $\Gamma_\tau$ is a geodesic flow of this metric, and is therefore a Riemannian geodesic flow. Fuller details may be inferred from a discussion of geodesibility in \cite{Recht}.

Each geodesic of $F_\tau$, for fixed $\tau$, has constant curvature in the ordinary sense, that is, constant geodesic curvature with respect to the Euclidean metric, but different geodesics have different curvatures. Each geodesic also has constant geodesic curvature with respect to the Poincar\'{e} metric, but now all geodesics have the same curvature.
It may be shown (again by making use of Zermelo navigation) that $\tau=0,\pm 1$ are the only cases where $F_\tau$ is of constant flag curvature.

\subsubsection*{Acknowledgements}
The first author is a Guest Professor at Ghent University:\ he
is grateful to the Department of Mathematics for its
hospitality.

We are most grateful to Juan-Carlos Alvarez Paiva for a number of  observations, both helpful and challenging, and in particular for bringing the problem of finding a Finsler metric for the horocycle path space to our attention; and to David Saunders for his encouragement.

This work is part of the {\sc irses} project {\sc geomech} (nr.\
246981) within the 7th European Community Framework Programme.

\end{document}